\documentclass[12pt]{article}
\newcommand{\bm}[1]{\mbox{\boldmath $#1$}}      

\begin{document}

\newtheorem{prop}{Proposition}
\newtheorem{cor}{Corollary}
\newtheorem{lem}{Lemma}
\newtheorem{thm}{Theorem}

\newcommand{\e}[0]{e}                   
\newcommand{\fb}[0]{\raisebox{.6ex}{\framebox[0.5em]{}}}
\newcommand{\itg}[0]{\mbox{\bf Z}}              
\newcommand{\cpx}[0]{\mbox{\bf C}}              



\title
{\bf
Quantizing 
the B\"acklund transformations of 
Painlev\'e equations and 
the quantum discrete Painlev\'e VI equation
}
\author{Koji HASEGAWA
\vspace{1ex}
\\
Mathematical Institute, Tohoku University\\
Sendai 980-8578 JAPAN
\thanks{e-mail address : kojihas@math.tohoku.ac.jp}}
%
\date{}
%
%
%
\maketitle
\begin{center}
{\sl
Dedicated to Professor Akihiro Tsuchiya 
 on the occasion of his retirement} 
\end{center}
\begin{center}
{\bf Abstract}
\end{center}
{
Based on
 the 
works by Kajiwara, Noumi and Yamada,
we propose a canonically quantized version 
of the rational Weyl group representation 
which originally arose as 
``symmetries" or the B\"acklund transformations in 
Painlev\'{e} equations. 
We thereby propose a quantization of 
discrete Painlev\'{e}  VI equation as a discrete Hamiltonian flow 
commuting with the action of 
 $W(D_4^{(1)})$. 
%
%
}

\section{\bf Introduction}\label{sec:intro}

Let $A=[a_{ij}]_{i,j=0}^l$ be a generalized Cartan matrix of affine
type and $W(A)$ be the corresponding Weyl group. We denote the
generators by $s_i$ $(i=0,\cdots, l)$. 
%
Let 
 ${\bm F}_{cl}:=\cpx(a_0,\cdots,a_l, f_0, \cdots, f_l)$ 
be the field of rational functions generated by commuting variables 
  $a_0,\cdots,a_l, f_0, \cdots, f_l$.
Let $u_{ij}$ be integers that satisfy 

(i)
 $u_{ij}=0$ if $i=j$ or $a_{ij}=0$, 

(ii)
 $u_{ij}:u_{ji}=-a_{ij}:a_{ji}$ \quad otherwise.


\begin{thm}[KNY, case $A_l^{(1)}$]
For $i,j=0,\cdots, l$ put
\begin{equation}
s_i(a_j):=a_ja_i^{-a_{ij}}, 
\quad
s_i(f_j):=f_j\left(\frac{a_i+f_i}{1+a_if_i}\right)^{u_{ij}}.
\label{eq:KNYaction}
\end{equation}
Then these formulas 
define a group homomorphism 
$
W(A_l^{(1)}) \rightarrow {\rm Aut}(\bm F_{cl})
$.
\end{thm}

This is the typical formula of the affine Weyl group symmetry 
or the B\"acklund transforation 
for difference Painlev\'{e} equation in its symmetric form:
the case $A_2^{(1)}$ 
gives the symmetry of the difference 
Painlev\'{e} IV equation\cite{KNY}. 
Moreover, this action is a Poisson map with respect to the bracket
\begin{equation}
\{f_i, f_j\}=u_{ij}f_if_j,
\quad 
\{a_i, a_j\}=\{a_i, f_j\}=0.
\label{eq:poissonbrackets}
\end{equation}
A naive expect 
 is that there exist a quantization of this
representation realized as adjoint actions of some suitable operators
(quantum Hamiltonian action).
One of the 
aim of this note 
 is to 
answer this problem. 
In the type $A$ case, we introduce the 
letters $F_0, \cdots, F_l$ subject to the quantized relation
 of (\ref{eq:poissonbrackets}),
$$
F_iF_{i+1}=q^{-1}F_{i+1}F_i, \quad F_iF_j-F_jF_i=0\; (j\not\equiv i\pm 1)
$$
as well as central letters $a_0, \cdots, a_l$. 
Let $\bm F$ be the skew field defined by these relations. 
We will construct the affine Weyl group action on $\bm F$ in the form
$$
s_i(\phi)=S_i\phi S_i^{-1}
$$
for any $\phi\in \bm F$. 
The ``Hamiltonian" $S_i$ is actually given by some infinite product 
which is rather familiar in $q$- analysis, 
despite that it involves non-commutative letters
(Section \ref{sect:typeA}, Theorem \ref{thm:Abraidrel}).
It is also shown that 
the construction works for other affine Weyl groups as well
(Section \ref{sect:typeBG}, Theorems \ref{thm:B2case} and \ref{thm:G2case}). 

 
\medskip

Rescent studies of Painlev\'{e} systems 
enabled us to understand 
their discrete symmetries
(B\"acklund transformations) 
and the (discrete) time evolution tranformation 
from the one, 
namely the affine Weyl group actions of the above type
\cite{NY}\cite{Sakai}. 
Based on this knowledge together with our quantization of 
the affine Weyl group action, 
we can quantize
 discrete (multiplicative) Painlev\'e type equations. 
In principle, 
if we choose some lattice direction in the affine Weyl group as the
 generator of a discrete time evolution, 
then this discrete dynamics commutes with the 
 simple reflections corresponding to the roots that are
perpendicular to the evolution direction in the lattice.  
We apply this idea to quantize the
 $q$- difference Painlev\'e III equation studied by Kajiwara and Kimura
 \cite{KK}
and also  
to quantize Jimbo-Sakai's 
$q$- difference Painlev\'e VI system
\cite{JS}:
See 
(\ref{eq:qPIII4a_i}), (\ref{eq:qPIII4F_i}) 
in Section \ref{sect:typeA}
and 
Theorems \ref{thm:affD5action}, \ref{thm:qPVI} for results. 

\section{Quantizing the Weyl group action : Type A case}
\label{sect:typeA}
Let $Q={\bm Z}\alpha_0 +\cdots +{\bm Z}\alpha_{l}$ be the root
lattice of type $A_l^{(1)}$ 
with simple roots $\alpha_0, \cdots, \alpha_l$ and
$\bm C[Q]=\bm C[e^{\alpha_0}, \cdots, e^{\alpha_l}]$ 
be its group algebra. 
The Weyl group action $s_i(\alpha_j)=\alpha_j-a_{ij}\alpha_i$ gives rise to 
the action on 
$\bm C[Q]$, 
and we can identify the previously used letter $a_j$ as
$e^{\alpha_j}$:
$$
s_i(a_j)=a_i^{-a_{ij}}a_j.
$$
Let 
$\bm K$ be the quotient field of the group algebra
${\bm C}[Q]$, namely $\bm K={\bm C}(a_0,\cdots,a_l)$. 

\medskip

{\bf Remark}\quad 
Let
 $\partial_0, \cdots, \partial_l$ 
be the ``dual" letters such that 
$$
[\partial_j, \alpha_k]:=\partial_j\alpha_k-\alpha_h\partial_j=a_{jk},
$$
then we have
$$
e^{\pi\sqrt{-1}\alpha_i\partial_i}
\cdot \alpha_j \cdot
e^{-\pi\sqrt{-1}\alpha_i\partial_i}
=s_i(\alpha_j).
$$
That is, the Weyl group action on $\bm K$ can be realized as adjoint actions.
This remark applies to the latter cases as well. \fb

\medskip

For type $A_{l}^{(1)}$ case $(l>1)$, we introduce the 
cannonically quantized letters $F_0, \cdots, F_l$
corresponding to (\ref{eq:poissonbrackets}):
\begin{equation}
F_iF_{{i+1} {\rm mod}\,l+1}=q^{-1}F_{{i+1} {\rm mod}\, l+1}F_i, \; 
  F_iF_j=F_jF_i \quad (i-j\not\equiv \pm 1).
\label{eq:quantizedAalg}
\end{equation}

(Here and in what follows we regard the subscripts as elements in
 $\bm Z/(l+1)\bm Z$.)

\smallskip
\noindent
Let 
$\bm K\langle F_0,\cdots,F_l\rangle$
be the $\bm K$- algebra generated by the above letters
(\ref{eq:quantizedAalg}). It can be shown in a standard way that 
this algebra is an Ore domain (cf. \cite{B}). 
Let 
$\bm F:=\bm K(F_0,\cdots,F_l)$ be the quotient skew field of
 $\bm K\langle F_0,\cdots,F_l\rangle$. 
The above 
relations 
(\ref{eq:quantizedAalg}) 
actually quantize the Poisson bracket 
(\ref{eq:poissonbrackets}). 
In fact, 
letting $q\rightarrow 1$ and think of the Poisson structure 
$$
\{\phi, \psi\}:=\lim_{q\rightarrow 1}\frac{1}{q-1}[\phi, \psi]
$$
on the commutative algebra
 $\bm F \; {\rm mod} (q-1)$. Then we have
$$
\{F_i, F_{i+1}\}\equiv 
\frac{1}{q-1}(F_iF_{i+1}-F_{i+1}F_i)
=
-F_iF_{i+1}
$$
according to the defining relation (\ref{eq:quantizedAalg}).

Note that
$a_i\in\bm K$ is central in $\bm F$.
Let us introduce the following multiplication operator
\begin{equation}
\Psi(z,F_i)=\Psi_q(z,F_i)
:=
\frac
{(qF_i,q)_\infty (F_i^{-1},q)_\infty}
{(zqF_i,q)_\infty (zF_i^{-1},q)_\infty}
\label{eq:Amultop}
\end{equation}
where 
$z$ and $q$ are central letters and
$(x,q)_\infty := \prod_{m=0}^\infty (1+xq^m)$. 
The right hand side of (\ref{eq:Amultop}) should be 
understood 
in the $q$- adic completion $\bm F((q))$ of $\bm F$. 
We put
$$
\rho_i:=e^{\frac{1}{2}\pi\sqrt{-1}\alpha_i\partial_i},
\quad 
S_i:=\Psi(z,F_i) \rho_i.
$$
Note that $\rho_i$ 
commutes with the variables $F_j,\; 0\leq j \leq l.$
We are interested in the adjoint action of $S_i$, 
$
Ad(S_i):\phi \in \bm F((q)) \mapsto S_i\phi S_i^{-1}\in\bm F((q)).
$

The statement of the following theorem essentially goes back to \cite{FV}.

\begin{thm}
We have
\begin{equation}
Ad(S_i)^2=id,
\label{eq:Areflid}
\end{equation}
\begin{equation}
S_iS_j=S_jS_i \quad (j\not\equiv i\pm 1),
\nonumber
\end{equation}
\begin{equation}
S_iS_{i+1}S_i=S_{i+1}S_iS_{i+1}
\label{eq:Abraidrel}
\end{equation}
where the index should read modulo $l+1$.
Hence $s_i \mapsto Ad(S_i)$ defines 
a group homomorphism
$$W(A_l^{(1)}) \rightarrow {\rm Aut}{\left(\bm F((q))\right)}.$$
\label{thm:Abraidrel}
\end{thm}


Let us calculate 
$
Ad(S_i)F_j=S_i F_{j} S_i^{-1}
$ 
first. We have
\begin{eqnarray}
Ad(S_i)(F_{i-1})
&=&
\frac{1+a_iF_i}{a_i+F_i}F_{i-1}
\nonumber\\
Ad(S_i)(F_{i+1})
&=&
F_{i+1}\frac{a_i+F_i}{1+a_iF_i}
\label{eq:quatizedAaction}
\\
Ad(S_i)(F_j)&=&F_j \quad (i-j\not\equiv\pm 1)
\nonumber
\end{eqnarray}
In fact these are the defining recurrence relation for the multiplication
operator 
(\ref{eq:Amultop}), that is, one can recover the formula
 of $\Psi(a_i, F_i)$ 
 from these modulo pseudo constants.
  
Now we can check ${\rm Ad}S_i^2=id$(\ref{eq:Areflid}) from these formula. 
Equivalently, 
since 
$$
S_i^2=\Psi(a_i, F_i)\Psi(a_i^{-1}, F_i)\rho_i^2,
$$ 
(\ref{eq:Areflid}) follows from the fact that
$
\Psi(a_i, F_i)\Psi(a_i^{-1}, F_i)
$
is a pseudo constant:
\begin{equation}
\Psi(a_i, F_i)\Psi(a_i^{-1}, F_i)
=
\Psi(a_i, qF_i)\Psi(a_i^{-1}, qF_i).
\label{eq:psipsi'ispseudoconst}
\end{equation}

As for (\ref{eq:Abraidrel}), we can do the similar computation to 
check it as
 the adjoint action 
on $\bm F((q))$, 
namely, compare the result when adjointly applied to generaters $F_i$.
However, 
(\ref{eq:Abraidrel}) is satisfied as an identity of elements in
$\bm F((q))$.
Actually (\ref{eq:Abraidrel}) 
follows from  
the dilogarithmic identity 
: 
suppose $F$ and $G$ satisfies $FG=qGF$, then we have 
$$
(G,q)_\infty(F,q)_\infty=(F,q)_\infty(GF,q)_\infty(G,q)_\infty.
$$
>From this we can 
show 
 (\cite{FV}, \cite{Ki})
\begin{equation}
\Psi(x, F_i)\Psi(xy, F_{i+1})\Psi(y, F_i)
=
\Psi(y, F_{i+1})\Psi(xy, F_i)\Psi(x, F_{i+1})
\label{eq:dilogid}
\end{equation}
where $x, y$ are central letters, 
which is equivalent to  
(\ref{eq:Abraidrel}). \fb

\bigskip

%
Introduce the diagram automorphism by
$$
\omega: 
a_i\mapsto a_{i+1}, \quad F_i \mapsto F_{{i+1}\;{\rm mod}\,l+1},
$$
then $\omega$ and 
$s_i:=Ad(S_i)$ generate the extended affine Weyl group $\tilde{W}(A_{l}^{(1)})$ 
acting on $\bm F((q))$.
As is well known, we have the commuting elements 
\begin{equation}
\left\{
\begin{array}{ccl}
T_1&:=& s_1s_2\cdots s_l\omega^{-1}
\\
&&\\
T_2&:=& s_2\cdots s_l\omega^{-1} s_1
\\
& \vdots &
\\
T_l&:=& s_l\omega^{-1}s_1\cdots s_{l-1}.
\end{array}
\right.
\label{eq:commeltsinAbraid}
\end{equation}
They are mutually conjugate. If we take $T_1$ as a discrete time evolution 
operator, then the group 
$
\langle s_0s_1s_0, s_2, s_3, \cdots, s_l \rangle \simeq W(A_{l-1}^{(1)})
$
commutes with the $T_1$ action. 
This gives the quantization of 
the ``q- difference" version of type A discrete system 
with Painlev\'e type symmetry which is extensively studied by 
Noumi and Yamada \cite{NY}\cite{NY-A}.

\medskip

{\bf Example}\quad 
Let $l=2$. 
Note that 
$a_0a_1a_2=:p$ is invariant under $\tilde{W}(A_2^{(1)})$, 
and the same holds for 
$F_0F_1F_2=:c$ since $c$ commutes with everything.
The action of $T_1=s_1s_2\omega^{-1}$ is given by 
\begin{equation}
T_1(a_0)=p^{-1}a_0, \; 
T_1(a_1)=pa_1, \; 
T_1(a_2)=a_2 
\label{eq:qPIII4a_i}
\end{equation}
and
\begin{equation}
F_0T_1(F_0)=c\frac{1+a_1F_1^{-1}}{1+a_1F_1},
\quad
T_1^{-1}(F_1)F_1=c\frac{1+a_0^{-1}F_0}{1+a_0^{-1}F_0^{-1}}
.
\label{eq:qPIII4F_i}
\end{equation}
This 
$T_1$ action commutes with $\langle s_0s_1s_0, s_2\rangle\simeq W(A_1^{(1)})$ 
and
gives a quantization of the $q{\rm P_{III}}$ system studied in \cite{KK} 
(where $p$ should be regarded as $q$).
%

In this $l=2$ case, the  Hamiltonian for the diagram automorphism $\omega$ can be 
found as follows, 
so that the above  $T_1$ flow is actually a discrete Hamiltonian flow. 
We put 
$$
\theta(X):=(X,q)_\infty(qX^{-1},q)_\infty
$$
and
$$
\Omega:=\left(
\theta(F_0^{-1}F_1)\theta(qF_1)^2\theta(F_2^{-1}F_0^{-1})
\right)^{-1}
\times
p^{-\partial'_1}\rho_1\rho_2,
$$
where 
$$
[\partial'_1, \alpha_0]=-1, 
[\partial'_1, \alpha_1]=1, 
[\partial'_1, \alpha_2]=0.
$$

\noindent
(Note : 
If we realize the $A_2^{(1)}$ root system in $\bm R^3\oplus\bm R\delta$ 
by 
$
\alpha_0=e_3-e_1+\delta, \alpha_1=e_1-e_2, \alpha_2=e_2-e_3,
$
where $e_1, e_2, e_3$ are the standard orthogonal basis of $\bm R^3$ and 
$\delta$ 
the canonical null root 
so that $p=e^\delta$,
then $\partial'_1$ stands for the derivation corresponding to $e_1$. )
Then we can easily check
$$
\Omega a_i \Omega^{-1}= a_{i+1 \,{\rm mod}\, 3}, \;
\Omega F_i \Omega^{-1}
=F_{i+1 \,{\rm mod}\, 3}
$$
and therefore
\begin{equation}
T_1={\rm Ad}(S_1S_2\Omega^{-1}).
\label{eq:Hamiltonian4qPIII}
\end{equation}
We have 
\begin{eqnarray*}
S_1S_2\Omega^{-1}
&=&
\Psi(a_1,F_1)\rho_1
\Psi(a_2,F_2)\rho_2
\times(p^{-\partial'_1}\rho_1\rho_2)^{-1}
\theta(F_0^{-1}F_1)\theta(qF_1)^2\theta(F_2^{-1}F_0^{-1})
\\
&=&
\Psi(a_1,F_1)
\Psi(a_1a_2,F_2)
\theta(F_0^{-1}F_1)\theta(qF_1)^2\theta(F_2^{-1}F_0^{-1})
\times p^{\partial'_1}.
\end{eqnarray*}
\quad\hfill\fb

\section{General case}
\label{sect:typeBG}
If the Dynkin diagram for the genralized Cartan matrix is simply laced, 
the construction in the last section applies to obtain the 
corresponding Weyl group action. 
As for the non-simply laced case, the construction 
can be reduced to the rank two cases : $B_2$ type and $G_2$ type 
(cf. \cite{NY-W}).

As before let $\bm K:=\bm C(a_1,a_2)$ be the 
quotient field of the group algebra $\bm C[Q]$, 
where $Q$ stands for the rank two root lattice in problem 
and we identify the letter
 $a_j$ 
with 
$e^{\alpha_j}\in \bm C[Q]$.  
We introduce 
 $\partial_j (j=1,2)$ such that $[\partial_j, \alpha_k]=a_{jk}$ and put
$
\rho_j:=e^{\frac{1}{2}\pi\sqrt{-1}\alpha_j\partial_j}. 
$
Then the Weyl group action $s_j$ on $\bm K$ 
is given by the adjoint action of $\rho_j$: 
$$
s_j(\alpha_k):=
\alpha_k-a_{jk}\alpha_j=
\rho_j\alpha_k\rho_j^{-1}, 
\quad
s_j(a_k)=a_j^{-a_{jk}}a_k=
\rho_ja_k\rho_j^{-1}.
$$
%

\underline{$B_2$ case.}

Let 
$
\left[\begin{array}{cc}
a_{11} & a_{12}\\a_{21} & a_{22}
\end{array}\right]
=
\left[\begin{array}{rr}
 2 & -1 \\ -2 & 2
\end{array}\right]
$
and define the skew field 
$\bm F=\bm K(F_1, F_2)$, where
$$
F_2F_1 =q^{2}F_1F_2.
$$

\underline{$G_2$ case.}

Let 
$
\left[\begin{array}{cc}
a_{11} & a_{12}\\a_{21} & a_{22}
\end{array}\right]
=
\left[\begin{array}{rr}
 2 & -1 \\ -3 & 2
\end{array}\right]
$
and define the skew field 
$\bm F=\bm K(F_1, F_2)$, where
$$
F_2F_1 =q^{3}F_1F_2.
$$




Using these, 
we have Hamiltonian Weyl group action on $\bm F$ in both cases :

\begin{thm}For type $B_2$ case, put 
\begin{equation}
S_1:=\Psi_q(a_1,F_1)\rho_1, 
\quad
S_2:=\Psi_{q^2}(a_2^2,F_2) \rho_2. 
\label{eq:S_i4Bcase}
\end{equation}
Then we have
$$Ad(S_j)^2=id, $$
\begin{equation}
S_1S_2S_1S_2=S_2S_1S_2S_1. 
\label{eq:B2case}
\end{equation}
\label{thm:B2case}
\end{thm}

\begin{thm}
For type $G_2$ case, put
\begin{equation}
S_1:=\Psi_q(a_1,F_1) \rho_1, 
\quad
S_2:=\Psi_{q^3}(a_2^3,F_2) \rho_2. 
\label{eq:S_i4Gcase}
\end{equation}
We have 
$$
Ad(S_j)^2=id, 
$$
\begin{equation}
S_1S_2S_1S_2S_1S_2=S_2S_1S_2S_1S_2S_1. 
\label{eq:G2case}
\end{equation}
\label{thm:G2case}
\end{thm}

For example, in the $B_2$ case we have
\begin{eqnarray*}
S_1 F_2  S_1^{-1}
&=&
\Psi_q(a_1,F_1) F_2 \Psi_q(a_1,F_1)^{-1}
\\&=&
F_2 \Psi_q(a_1,q^{-2}F_1) \Psi_q(a_1,F_1)^{-1}
\\&=&
F_2 
\frac
{(F_1q^{-1},q)_\infty(F_1^{-1}q^2,q)_\infty}
{(a_1F_1q^{-1},q)_\infty(a_1F_1^{-1}q^2,q)_\infty}
\frac
{(a_1F_1q,q)_\infty(a_1F_1^{-1},q)_\infty}
{(F_1q,q)_\infty(F_1^{-1},q)_\infty}
\\&=&
F_2 
\frac
{(1+F_1q^{-1})(1+F_1)}
{(1+a_1F_1q^{-1})(1+a_1F_1)}
\frac
{(1+a_1F_1^{-1})(1+a_1F_1^{-1}q)}
{(1+F_1^{-1})(1+F_1^{-1}q)}
\\&=&
F_2 
\frac
{(a_1+F_1)(a_1+F_1q^{-1})}
{(1+a_1F_1)(1+a_1F_1q^{-1})}.
\end{eqnarray*}
The property $Ad(S_j)^2=id$ can be checked by continuing this 
computation, 
or we can conclude it immediately from the pseudoconstant property 
(\ref{eq:psipsi'ispseudoconst}).

In principle we can calculate 
the expressions for 
$Ad(S_1S_2S_1S_2)F_j$ and $Ad(S_2S_1S_2S_1)F_j$ $(j=1,2)$ 
also 
to check 
(\ref{eq:B2case}) at the 
adjoint level, though it is a quite lengthy way.
In fact 
 (\ref{eq:B2case}) holds as an identity in $\bm F((q)).$ 
Let us introduce a 
square root of $F_2$, 
namely let 
$\sqrt{-F_2}$
 be the letter satisfying
$$
\sqrt{-F_2}^2=-F_2,\quad F_1\sqrt{-F_2}=q^{-1}\sqrt{-F_2}F_1.
$$
Then we have
$$
\Psi_{q^2}(a_2^2, F_2)
=
\Psi_{q}({a_2}, \sqrt{-F_2})\Psi_{q}({a_2}, -\sqrt{-F_2}),
$$
so that 
(\ref{eq:B2case}) 
 can be reduced to the type $A$ identity 
(\ref{eq:dilogid})
. 
For short, 
we write 
$a_1=:a, a_2=:b$, 
$$
\Psi_q(a, F_1)=:\Psi_1^{a}, \;
\Psi_{q^2}(b^2, F_2)=:\Psi_2^b, \;
\Psi_q(b, \pm\sqrt{-F_2})=:\Psi_\pm^b.
$$
Then $S_1S_2S_1S_2=S_2S_1S_2S_1$ is equivalent to
$
\Psi_1^a\Psi_2^{ab}\Psi_1^{ab^2}\Psi_2^b
=
\Psi_2^b\Psi_1^{ab^2}\Psi_2^{ab}\Psi_1^{a},
$
or
$$
\Psi_1^a\Psi_+^{ab}\Psi_-^{ab}\Psi_1^{ab^2}\Psi_+^b\Psi_-^b
=
\Psi_+^b\Psi_-^b\Psi_1^{ab^2}\Psi_+^{ab}\Psi_-^{ab}\Psi_1^{a}.
$$
This can be verified as follows, 
which uses
(\ref{eq:dilogid}):
$\Psi_\pm^x\Psi_1^{xy}\Psi_\pm^y
=\Psi_1^y\Psi_\pm^{xy}\Psi_1^x
$ 
as well as
 $\Psi_+^x\Psi_-^y=\Psi_-^y\Psi_+^x$ 
at the underlined places. 
\begin{eqnarray*}
{\rm LHS}&=&
\Psi_1^a\Psi_+^{ab}
\underline{\Psi_-^{ab}\Psi_1^{ab^2}\Psi_-^b}
\Psi_+^b
\\
&=&
\underline{\Psi_1^a\Psi_+^{ab}\Psi_1^{b}}
\Psi_-^{ab^2}\Psi_1^{ab}\Psi_+^b
\\
&=&
\Psi_+^b\Psi_1^{ab}
\underline{\Psi_+^{a}\Psi_-^{ab^2}}
\Psi_1^{ab}\Psi_+^b
=
\Psi_+^b\Psi_1^{ab}\Psi_-^{ab^2}
\underline{
\Psi_+^{a}\Psi_1^{ab}\Psi_+^b
}
\\
&=&
\Psi_+^b
\underline{\Psi_1^{ab}\Psi_-^{ab^2}\Psi_1^{b}}
\Psi_+^{ab}\Psi_1^a
\\
&=&
\Psi_+^b
\Psi_-^{b}\Psi_1^{ab^2}\Psi_-^{ab}
\Psi_+^{ab}\Psi_1^a
=
{\rm RHS}.
\end{eqnarray*}

\medskip

Proof for the $G_2$ case 
can be quite similary done as in the $B_2$ case : we use
the cubic root
 $\zeta\neq 1$ of unity and $\sqrt[3]{F_2}$ of $F_2$ 
that satisfies $F_1\sqrt[3]{F_2}=q^{-1}\sqrt[3]{F_2}F_1$. 
We have
$$
\Psi_{q^3}(a_2^3,F_2)=
\Psi_{q}\left({a_2}, \sqrt[3]{F_2}\right) 
\Psi_{q}\left({a_2}, \zeta\sqrt[3]{F_2}\right) 
\Psi_{q}\left({a_2}, \zeta^{-1}\sqrt[3]{F_2}\right).
$$
As before, write $a_1=:a, a_2=:b$, 
$
\Psi_q(a, F_1)=:\Psi_1^{a}, 
\Psi_{q^3}(a^3, F_2)=:\Psi_2^b
$
and also
$$
\Psi_q(b, \sqrt[3]{F_2})=:\Psi_0^b, \quad
\Psi_q(b, \zeta^{\pm 1}\sqrt[3]{F_2})=:\Psi_\pm^b 
$$
for short. 
Then 
(\ref{eq:G2case})
is equivalent to
\begin{equation}
\Psi_1^a\Psi_2^{ab}\Psi_1^{a^2b^3}
\Psi_2^{ab^2}\Psi_1^{ab^3}\Psi_2^b
=
\Psi_2^b\Psi_1^{ab^3}\Psi_2^{ab^2}
\Psi_1^{a^2b^3}\Psi_2^{ab}\Psi_1^a.
\label{eq:restateG2}
\end{equation}
This time 
(\ref{eq:dilogid}) means  
$
\Psi_k^x\Psi_1^{xy}\Psi_k^y=\Psi_1^y\Psi_k^{xy}\Psi_1^x
$
for $k=0,\pm$ and
$
\Psi_0^x, \Psi_+^y, \Psi_-^z
$
are commuting for any central $x, y, z$.
We have
\begin{eqnarray*}
{\rm LHS \; of \; (\ref{eq:restateG2})}
&=&
\Psi_1^a
\Psi_0^{ab}\Psi_+^{ab}
\underline{\Psi_-^{ab}
\Psi_1^{a^2b^3}}
\underbrace{\Psi_0^{ab^2}}\Psi_+^{ab^2}\underline{\Psi_-^{ab^2}}
\underbrace{\Psi_1^{ab^3}
\Psi_0^b}
\Psi_+^b\Psi_-^b
\\&=&
\Psi_1^a
\Psi_0^{ab}\Psi_+^{ab}
\Psi_1^{ab^2}
\Psi_-^{a^2b^3}
\underline{
\Psi_1^{ab}
\Psi_+^{ab^2}
\Psi_1^b
}
\Psi_0^{ab^3}
\Psi_1^{ab^2}
\Psi_+^b\Psi_-^b
\\&=&
\Psi_1^a
\Psi_0^{ab}
\underline{
\Psi_+^{ab}
\Psi_1^{ab^2}
}
\Psi_-^{a^2b^3}
\underline{
\Psi_+^{b}
}
\Psi_1^{ab^2}
\underbrace{
\Psi_+^{ab}
}
\Psi_0^{ab^3}
\underbrace{
\Psi_1^{ab^2}
\Psi_+^b
}
\Psi_-^b
\\&=&
\underline{
\Psi_1^a
\Psi_0^{ab}
\Psi_1^{b}
}
\Psi_+^{ab^2}
\Psi_1^{ab}
\Psi_-^{a^2b^3}
\underbrace{
\Psi_1^{ab^2}
\Psi_0^{ab^3}
\Psi_1^{b}
}
\Psi_+^{ab^2}
\Psi_1^{ab}
\Psi_-^b
\\&=&
\Psi_0^b
\Psi_1^{ab}
\underline{
\Psi_0^{a}
}
\Psi_+^{ab^2}
\underline{
\Psi_1^{ab}
}
\Psi_-^{a^2b^3}
\underline{
\Psi_0^{b}
}
\Psi_1^{ab^3}
\Psi_0^{ab^2}
\Psi_+^{ab^2}
\Psi_1^{ab}
\Psi_-^b
\\&=&
\Psi_0^b
\underline{
\Psi_1^{ab}
\Psi_+^{ab^2}
\Psi_1^{b}
}
\Psi_0^{ab}
\underbrace{
\Psi_1^{a}
\Psi_-^{a^2b^3}
\Psi_1^{ab^3}
}
\Psi_0^{ab^2}
\Psi_+^{ab^2}
\Psi_1^{ab}
\Psi_-^b
\\&=&
\Psi_0^b
\Psi_+^{b}
\Psi_1^{ab^2}
\Psi_+^{ab}
\underline{
\Psi_0^{ab}
}
\Psi_-^{ab^3}
\underline{
\Psi_1^{a^2b^3}
}
\underbrace{
\Psi_-^{a}
}
\underline{
\Psi_0^{ab^2}
}
\Psi_+^{ab^2}
\underbrace{
\Psi_1^{ab}
\Psi_-^b
}
\\&=&
\Psi_0^b
\Psi_+^{b}
\Psi_1^{ab^2}
\Psi_+^{ab}
\Psi_-^{ab^3}
\underline{
\Psi_0^{ab}
}
\underline{
\Psi_1^{a^2b^3}
}
\underline{
\Psi_0^{ab^2}
}
\Psi_+^{ab^2}
\underbrace{
\Psi_-^{a}
}
\underbrace{
\Psi_1^{ab}
\Psi_-^b
}
\\&=&
\Psi_0^b
\Psi_+^{b}
\Psi_1^{ab^2}
\Psi_+^{ab}
\Psi_-^{ab^3}
\Psi_1^{ab^2}
\Psi_0^{a^2b^3}
\underline{
\Psi_1^{ab}
\Psi_+^{ab^2}
\Psi_1^{b}
}
\Psi_-^{ab}
\Psi_1^a
\\&=&
\Psi_0^b
\Psi_+^{b}
\Psi_1^{ab^2}
\underline{
\Psi_+^{ab}
}
\Psi_-^{ab^3}
\underline{
\Psi_1^{ab^2}
}
\Psi_0^{a^2b^3}
\underline{
\Psi_+^{b}
}
\Psi_1^{ab^2}
\Psi_+^{ab}
\Psi_-^{ab}
\Psi_1^a
\\&=&
\Psi_0^b
\Psi_+^{b}
\underline{
\Psi_1^{ab^2}
\Psi_-^{ab^3}
\Psi_1^{b}
}
\Psi_+^{ab^2}
\underbrace{
\Psi_1^{ab}
\Psi_0^{a^2b^3}
\Psi_1^{ab^2}
}
\Psi_+^{ab}
\Psi_-^{ab}
\Psi_1^a
\\&=&
\Psi_0^b
\Psi_+^{b}
{
\Psi_-^{b}
\Psi_1^{ab^3}
\Psi_-^{ab^2}
}
\Psi_+^{ab^2}
{
\Psi_0^{ab^2}
\Psi_1^{a^2b^3}
\Psi_0^{ab}
}
\Psi_+^{ab}
\Psi_-^{ab}
\Psi_1^a
\\&=&
\Psi_2^b\Psi_1^{ab^3}\Psi_2^{ab^2}
\Psi_1^{a^2b^3}\Psi_2^{ab}\Psi_1^a.
=
{\rm RHS \; of \; (\ref{eq:restateG2})}
\quad\qquad \hfill\fb
\end{eqnarray*}

\medskip



\section{Quantizing the discrete Painlev\'e equation}

There is a discretization of the Painlev\'e VI equation proposed by
Jimbo and Sakai \cite{JS} 
and later it was reformulated under the affine Weyl group symmetry 
of type $D_5^{(1)}$ : \cite{Sakai}, \cite{TM}.
Using the ideas in the previous sections, 
here we propose its quantum (non-commutative) 
version in a quite straightforward way. 
%
Let us introduce the Dynkin diagram of type $D_5^{(1)}$ and its 
numbering:
$$
\begin{array}{cccccccc}
 0\quad  &   & \quad 5\\
\backslash&   & \slash  \\
 \quad 2 & - & 3 \quad    \\
  \slash  &   & \backslash  \\
1 \quad  &   & \quad   4
\end{array}
$$
We denote the corresponding generalized Cartan matrix by 
$[a_{ij}]_{i,j=0}^{5}$ 
and the simple roots by $\{\alpha_i\}_0^5$.
Let 
 $q$ be a formal central letter (or a complex parameter, $|q|<1$). 
We introduce the 
field of rationals $\bm K=\bm C(a_0, \cdots, a_5)$, 
where $a_i=e^{\alpha_i}$ as in the previous sections.
%
%

The Weyl group $W=W(D_5^{(1)})$ acts on
 $\bm K$ 
by
$$
s_i(a_j)=a_i^{-a_{ij}}a_j.
$$
One checks that $a_0a_1a_2^2a_3^2a_4a_5=:p$ is invariant under the 
action of $W$. Moreover, this action can be extended 
by the diagram automorphisms
$
\sigma_{01}, \sigma_{45}, \tau: 
$
\begin{eqnarray}
\sigma_{01}:
a_0\leftrightarrow a_1^{-1}, a_j\mapsto a_j^{-1}\; (j\neq 0,1),
&&
\sigma_{45}:
a_4\leftrightarrow a_5^{-1}, a_j\mapsto a_j^{-1}\; (j\neq 4, 5),
\nonumber
\\
\tau : a_j &\leftrightarrow& a_{5-j}^{-1} \; (j=0,\cdots, 5).
\label{eq:diagautom}
\end{eqnarray}
We denote the extended Weyl group by
 $\tilde{W}:=\langle W, \sigma_{01}, \sigma_{45}, \tau\rangle$.

\bigskip

Let 
$\bm F=\bm K(F, G)$ be the 
skew field, 
where 
\begin{equation}
FG=qGF.
\label{eq:defP6field}
\end{equation}
The 
 action of $\tilde{W}$
on 
$\bm K$
can be extended to 
 $\bm F=\bm K(F, G)$ 
as follows.


\begin{thm}
We can extend the 
automorphisms
 $s_j$ ($j=0, \cdots, 5$) of  $\bm K$ 
as algebra automorphisms 
of $\bm F$  
by putting
$$
s_2(F):= F \frac{a_0a_1^{-1}G+a_2^2}{a_0a_1^{-1}a_2^2G+1}, \quad
s_j(F):= F \; (j\neq 2)
$$
and
$$
s_3(G):=   \frac{a_3^2a_4a_5^{-1}F+1}{a_4a_5^{-1}F+a_3^2} G, \quad
s_j(G):= G \; (j\neq 3).
$$
They give rise to 
a homomorphism
 $W \rightarrow {\rm Aut}_{\rm{ skew\; field}}(\bm F).$

Moreover, the action of the diagram automorphisms 
$
\sigma_{01}, \sigma_{45}, \tau 
$
on $\bm K$
can be extended as involutive antiautomorphisms on $\bm F$ by
$$
\sigma_{01} : F \mapsto q^{-1}F^{-1}, \; G \mapsto G
$$
$$
\sigma_{45} : F\mapsto F, \; G \mapsto q^{-1}G^{-1}
$$
$$
\tau : F \mapsto G, \; G \mapsto F
$$
so that
we have a homomorphism 
$\tilde W \rightarrow {\rm Aut}_{\bm C-lin}(\bm F)$. \quad\fb

\label{thm:affD5action}
\end{thm}

\bigskip

The proof is straightforward. 
Because of the noncommutativity, 
it seems inevitable to define 
$\sigma_{01}, \sigma_{45}, \tau$ 
 actions on $\bm F$
as antiautomorphisms. 
For example, if we want to extend 
the action of $\tau$ on $\bm K$ to $\bm F$ 
by
$\tau : F\leftrightarrow G$ (cf. \cite{TM}), 
this cannot be compatible with the 
relation $FG=qGF$ if we insist $\tau$ 
to be an automorphism:
$\tau(F)\tau(G)\neq q\tau(G)\tau(F).$  

Note that 
$
\sigma:=\sigma_{01}\sigma_{45}
$ 
is an automorphism of $\bm F$ satisfying 
$
\sigma s_j=s_{\sigma(j)}\tau
$, where 
$$
(\sigma(0),\sigma(1),\sigma(2),\sigma(3),\sigma(4),\sigma(5))
= (1,0,2,3,5,4).
$$
We are interested in the action of
$$
T_3=s_2s_1s_0s_2\sigma_{01}s_3s_4s_5s_3\sigma_{45}
=
s_2s_1s_0s_2s_3s_4s_5s_3\sigma 
\in \tilde W,
$$
since 
in the commutative case this recovers the discrete Painlev\'e VI 
system \cite{Sakai}.
Put $t:=a_3^2a_4a_5,$
then $t$ is invariant under 
$\langle s_0, s_1, s_2s_3s_2=s_3s_2s_3, s_4, s_5
\rangle
= W^{T_3}
$.


\bigskip

\begin{thm}
We have
$$
%
(T_3(a_0), T_3(a_1), T_3(a_2), T_3(a_3), T_3(a_4), T_3(a_5))
=
(a_0, a_1, pa_2, p^{-1}a_3, a_4, a_5), 
\;
T_3(t)=p^{-1}t,
$$
\begin{equation}
T_3(F)
=
q^{-1}p^2t^{-2}
\frac
{\;G+tp^{-1}a_1^2\;}
{\;G+t^{-1}pa_0^2\;}
\frac
{\;G+tp^{-1}a_1^{-2}\;}
{\;G+t^{-1}pa_0^{-2}\;}
F^{-1}
\label{eq:Fbar}
\end{equation}
\begin{equation}
T_3^{-1}(G)
=
q^{-1}t^{-2}G^{-1}
\frac{\;F+ta_4^2\;}{\;F+t^{-1}a_5^{2}\;}
\frac{\;F+ta_4^{-2}\;}{\;F+t^{-1}a_5^{-2}\;}
.
\label{eq:Gbar}
\end{equation}
This $T_3$ flow 
allows the symmetry of $W(D_4^{(1)})$, 
namely $T_3$ commutes with the subgroup 
$\langle s_0, s_1, s_2s_3s_2=s_3s_2s_3, s_4, s_5, 
\sigma_{01}, \tau 
\rangle 
\simeq 
\langle W(D_4^{(1)}), \sigma_{01}, \tau 
\rangle$
 of $\tilde W$. 
\qquad\fb
\label{thm:qPVI}
\end{thm}


Put $Z:=tF, Y:=\frac{t}{p}G, T:=\left(\frac{t}{p}\right)^2$ and 
let us use the notation $T_3(X)=:\bar{X}$.
We have $ZY=qYZ$, $\bar{T}=p^{-2}T$ and 
the formula (\ref{eq:Fbar}), (\ref{eq:Gbar}) can be rewritten as
follows:
\begin{equation}
\bar{Z}Z
=\frac{p}{q}
\frac{\,Y+Ta_1^2\,}{\,Y+a_0^2\,}
\frac{\,Y+Ta_1^{-2}\,}{\,Y+a_0^{-2}\,}
,
\qquad
\bar{Y}Y
=\frac{1}{pq}
\frac{\,\bar{Z}+Ta_4^2\,}{\,\bar{Z}+a_5^2\,}
\frac{\,\bar{Z}+Ta_4^{-2}\,}{\,\bar{Z}+a_5^{-2}\,}
%
.
\end{equation}
This system should be regarded as 
a quantization of the discrete Painlev\'e VI equation.

Likewise in the previous sections, we have the Hamiltonians for the 
$W(D_5^{(1)})$-
 action written in terms of infinite product $\Psi$. 
Put 
$$
S_2=\Psi(a_2^2, a_0a_1^{-1}G)\rho_2, \;
S_3=\Psi(a_3^2, a_5a_4^{-1}F)\rho_3, \;
S_j=\rho_j
\; (j\neq 2,3),
$$
where $[\partial_j, a_k]=a_{jk}a_k$
and 
$\rho_j=e^{\frac{\pi}{2}\sqrt{-1}\alpha_j\partial_j}.$
Then we have
$
s_j={\rm Ad}(S_j)
$
for $j=0, \cdots, 5.$

We can also find a Hamiltonian for the diagram automorphism
 $
 \sigma=\sigma_{01}\sigma_{45}.
 $ 
Let us introduce the letters $\partial'_j$ by the relation 
$[\partial'_j, \alpha_k]=\delta_{jk}$
and put
\begin{eqnarray*}
\Sigma
&:=&
%
\theta(qFG)\theta(G^{-1}F)\theta(qF)^4
\\&&
\times
e^{\frac{\pi}{2}\sqrt{-1}(\alpha_0+\alpha_1)(\partial'_0+\partial'_1)}
e^{{\pi}\sqrt{-1}\alpha_2\partial'_2}
e^{{\pi}\sqrt{-1}\alpha_3\partial'_3}
e^{\frac{\pi}{2}\sqrt{-1}(\alpha_4+\alpha_5)(\partial'_4+\partial'_5)}.
\end{eqnarray*}
Then we have 
$
\sigma=Ad(\Sigma), 
$
that is, 
$\Sigma F \Sigma^{-1}= F^{-1}$, 
$\Sigma G \Sigma^{-1}= G^{-1}$ 
and
$\Sigma a_j \Sigma^{-1}= a_{\sigma(j)}^{-1}$ 
hold.
Thus we have
\begin{thm}
The quantum discrete Painlev\'e VI equation is a discrete Hamiltonian flow, 
$$
T_3= Ad
(S_2S_1S_0S_2S_3S_4S_5S_3\Sigma
).
$$
\quad\hfill\fb
\label{thm:qP6isHam}
\end{thm}
Explicitly, we have
\begin{eqnarray*}
&&
S_2S_1S_0S_2S_3S_4S_5S_3\Sigma 
\\
&=&
\Psi(a_2^2,a_0a_1^{-1}G)\rho_2
\rho_1\rho_0
\Psi(a_2^2,a_0a_1^{-1}G)\rho_2
\Psi(a_3^2,a_4^{-1}a_5F)\rho_3
\rho_4\rho_5
\Psi(a_3^2,a_4^{-1}a_5F)\rho_3
\Sigma
\\
&=&
\Psi(a_2^2,a_0a_1^{-1}G)
\Psi((a_0a_1a_2)^2,a_0^{-1}a_1G)
\Psi(p^2(a_3a_4a_5)^{-2},a_4a_5^{-1}F)
\Psi(p^2a_2^2a_3^{-2},a_4a_5^{-1}F)
\\&&
\times
\theta(qFG)\theta(G^{-1}F)\theta(qF)^4
p^\partial, 
\end{eqnarray*}
where
$
\partial=
\frac{1}{2}(\partial'_3-\partial'_2)
$
so that 
$
p^\partial t p^{-\partial}=tp^{-1}.
$
Note that we can realize the $D_5^{(1)}$ root lattice 
in 
$\bm R\delta\oplus\bm R^5
=\bm R\delta\oplus\bm Re_1\oplus\cdots\oplus\bm Re_5$ 
by
$$
\alpha_0=\delta-e_1-e_2, \;
\alpha_1=e_1-e_2, \;
\alpha_2=e_2-e_3, \;
\alpha_3=e_3-e_4, \;
\alpha_4=e_4-e_5, \;
\alpha_5=e_4+e_5, 
$$
where $e_j$ are regarded as the orthonormal basis. 
Then we have 
$t=e^{e_3},$ 
$\partial=-\partial/\partial e_3.$
The root subsystem pependicular to $e_3$ is 
generated by 
$\alpha_0, \alpha_1, \alpha_2+\alpha_3=e_2-e_4, \alpha_4, \alpha_5$ 
and isomorphic to $D_4^{(1)}.$

\bigskip

As for the classical (commutative) case, 
the discrete Painlev\'e 
system allows rather simple,  
so-called ``seed" solutions, 
from which one can construct rich explicit solutions via
 B\"acklund transformations. 
In our quantized system, 
it 
seems however that 
such seed solutions are difficult to find 
because of the noncommutativity.
For such issues as well as the consideration of the continuous limit, 
we hope to discuss elsewhere.


%
%



\bigskip

{\bf Acknowledgement.} 
Part of 
this work is based on a talk at Newton institute, Cambridge
(EuroConference ``Application of the Macdonald polynomials'',
16 -21 April 2001) and thanks are due to the organizers, 
including Professor Noumi, of the conference.
He also express a sincere gratitude for Gen Kuroki, 
Tetsuya Kikuchi and Hajime Nagoya for valuable discussions and 
informations. 

This work has been supported by the grants-in-aid for scientific research, 
Japan Society for the Promotion of Science, 
no.12640005, no. 16540182.  
%

\end{document}